\documentclass[12pt,letterpaper]{amsart}
\usepackage{palatino, euler, epic,eepic, amssymb, xypic, floatflt,   
  microtype}
\usepackage{hyperref}
\usepackage{amsmath,amsthm,amssymb,graphicx}

\input xy
\xyoption{all}
 \newlength{\baseunit}               % the basic unit length
 \newcount{\numlines}                % depth of picture (in number of lines)
\usepackage{pdfsync,url, mathrsfs}

\numberwithin{equation}{section}

\newcommand{\tpoint}[1]{\vspace{2mm}\par \noindent \refstepcounter{equation}{\theequation.}
  {\bf #1. ---} }

\newcommand{\bpoint}[1]{\vspace{2mm}\par \noindent \refstepcounter{equation}{\theequation.}
  {\bf #1.} }

\newcommand{\cut}[1]{}
\newcommand\hidden[1]{}

\title{Expanders from Markov bases}

\date\today
\author{Alexander Engstr\"om}
\address{Department of Mathematics, Aalto University, Helsinki, Finland}
\email{alexander.engstrom@aalto.fi}

\begin{document}
\begin{abstract}
Diaconis and Sturmfels introduced an influential method to construct Markov chains using commutative algebra. One major point of their method is that infinite families of graphs are simultaneously proved to be connected by a single algebraic calculation. For large state spaces in the infinite families these Markov chains are not rapidly mixing and only ad hoc methods have been available to improve their mixing times. We provide a method to get rapid mixing by constructing expanders for the Diaconis-Sturmfels type Markov chains.
\end{abstract}

\maketitle

\setcounter{tocdepth}{1} % this just includes subsections

% \tableofcontents

\section{The model problem}

In this text we discuss how to build expanders for a certain class of model problems. These problems were originally solved using the Markov bases designed by Diaconis and Sturmfels \cite{ds}, and they were not rapidly mixing in this setting. The book by Drton, Sturmfels and Sullivant \cite{dss} is an excellent introduction to algebraic statistics.

For basic facts about expanders and mixing time, we refer to the survey by Hoory, Linial and Widgerson \cite{hlw}. 
The expanders in this paper are built using the zig-zag product introduced by Reingold, Vadhan and Widgerson \cite{rvw}.
Mixing time and expansion are properties of sequences of larger and larger graphs, not of single graphs. We are interested in graphs from Markov bases and they come naturally in sequences for fixed design matrices. In this paper we don't study sequences of graphs from different design matrices; allowing for any design matrices essentially doesn't restrict the sequence of graphs and it's hard to make useful statements in that general setting.

\bpoint{The model problem} We are provided the following data:

{\em (i)} A finite set $S$ of rational numbers in $[ 0,1 ) ^d.$

{\em (ii)} A full dimensional polytope $P$ in $\mathbf{R}^d.$

For every positive integer $m$ we want to construct a connected graph $G_m$ whose vertex set is
\[
(S+\mathbf{Z}^d) \cap mP
\]
and it should be as fast as possible to random walk on $G_m.$ The main result of this paper is an explicit construction of graphs $G_m$ that are expanders by a straight-forward application of the zig-zag product. 

\bpoint{Markov bases are not rapidly mixing for the model problem} \label{nm}
The vertices of $G_m$ in our model problem have coordinates and we can assign a length to each edge of $G_m$ from the Eulidean metric. We consider a wider context than of Markov bases: Assume that the edge lengths are bounded from above by some $\ell$ for all $G_m,$ even as $m\rightarrow \infty.$ Fix an hyperplane $H$ that divides $P$ into two equal volumes, such that the (codimension one) volume of $P \cap H$ is minimal among those hyperplanes. Removing all vertices of $(S+\mathbf{Z}^d) \cap mP$ that are within distance $\ell$ from $mH$ will cut the graph $G_m$ into two disjoint subgraphs of similar orders. As $m\rightarrow \infty $ the proportion of vertices removed from $G_m$ to cut it into two pieces tends to 0, and this proves via the Cheeger inequality that we don't have expansion. In \ref{algi} we make some remarks on this in practice.

\bpoint{The model problem in algebraic statistics} A typical problem in algebraic statistics is to random walk on all $k \times k$ contingency tables with fixed row and columns sums $m.$ This is encoded as $A\mathbf{x}=m\mathbf{b}$ where $A$ is the {\em design matrix}. Both it and $\mathbf{b}$ are kept fixed while $m$ is allowed to grow. The non-negative integer vectors $\mathbf{x}$ encode the contingency tables. The non-negativity condition provides linear inequalities for a usually not full-dimensional polytope, but after restricting to the correct ambient space we have an instance of our model problem. For different $m$ we might get different polytopes $P,$ but one can deduce that they fall into a finite number of versions of our model problem.

\bpoint{Relation to other work} Windisch have generalised and provided detailed analysis of several aspects of our results, and provided an independent construction of expanders for contingency tables in his preprint \cite{w}.

\section{The main construction}

\begin{figure}
\includegraphics[width=12cm]{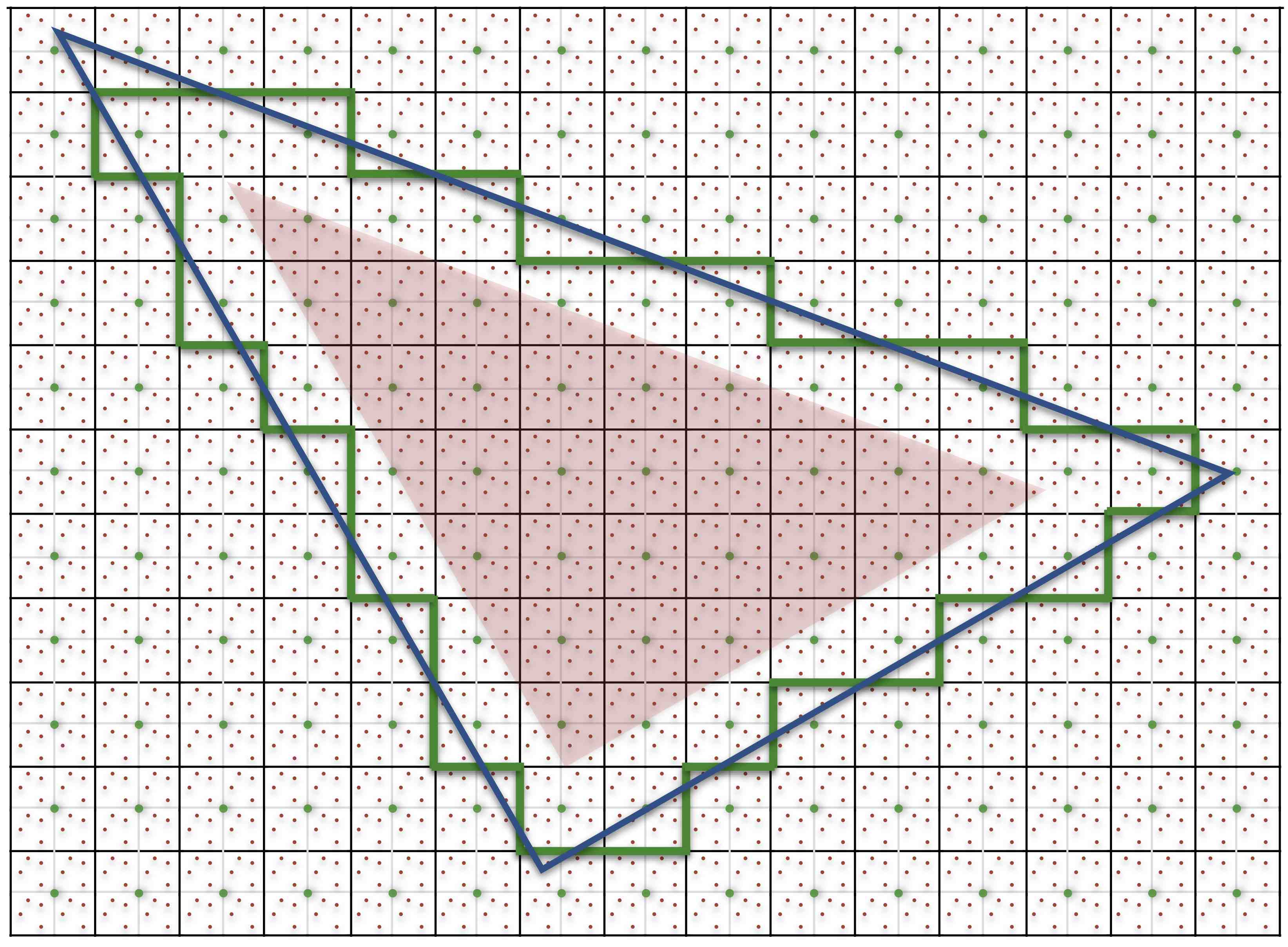}
\caption{Illustration of the main construction} \label{fig:B}
\end{figure}

\begin{figure}
\includegraphics[width=8cm]{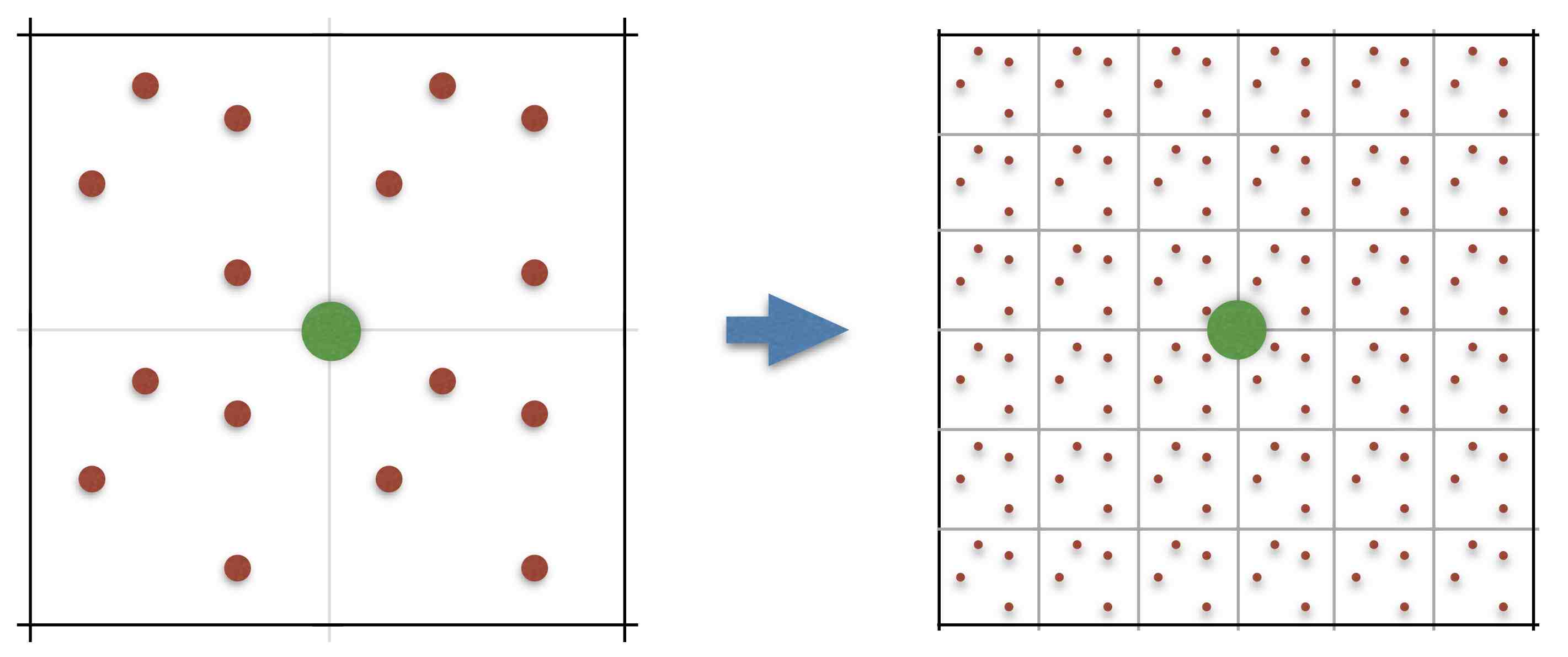}
\caption{Illustration of changing from $m=2$ to $m=6$ in the main construction.} \label{fig:A}
\end{figure}

We want to construct expanders for our model problem. In the construction we will refer to the example illustrated in Figures \ref{fig:B} and \ref{fig:A}.

\bpoint{Input data} We are provided the following data:

{\em (i)} A non-empty finite set $S=\{s_0,s_1, \ldots, s_{n_S-1}\}$ of rational numbers in $[ 0,1 ) ^d.$

{\em (ii)} A full dimensional polytope $P$ in $\mathbf{R}^d$ defined by rational linear inequalities. This is the pale red filled triangle in Figure~\ref{fig:B}.

For every positive integer $m$ we construct a graph $G_m.$ The first part of the construction is independent of $m$ and only calculated once -- this is similar to the spirit of a Markov basis.

\bpoint{The part independent of $m$}

{\em (iii)} Push each of the hyperplanes defining $P$ outwards by at least $\sqrt{d/2}$ to get a polytope $Q$ defined by rational linear inequalities. This is the blue triangle in Figure~\ref{fig:B}.

{\em (iv)} Define a connected $d_H$-regular graph $H$ whose $n_H$ vertices are in bijections to the set $ \{ z \in \mathbf{Z}^d \mid z + (1/2,\ldots,1/2) \in Q  \}.$ There are several properties of $H$ that would be beneficial but aren't necessary: Preferably the neighbours of a vertex should be possible to calculate fast at any vertex, instead of keeping $H$ in memory. It is also preferable if the vertex degree is small. The brute force way to construct $H$ is as a complete graph. The recommended and more gentle way is to use Markov bases, and if necessary add multiple self-loops at vertices to make it regular. The green big dots in Figure~\ref{fig:B} are the vertices of $H.$ 

{\em (v)} We denote the absolute value of the second largest eigenvalue of $H$ by $\lambda_H.$ It is advisable to have a reasonable upper bound for $\lambda_H$ at this point: $\lambda_H+2n_H^{-1/2}$ should be smaller than one for this construction to achieve expansion as $m \rightarrow \infty.$ In practice $\lambda_H$ is usually fine when $H$ is small and derived from a Markov basis. One can always make $H$ denser and push $\lambda_H$ towards $2n_H^{-1/2}.$ 

\bpoint{The part dependent of $m$}

{\em (vi)} Take a regular graph (expander) $E$ with minimal (absolute) second eigenvalue $\lambda_E$ provided that it should have $n_E=m^d$ vertices and vertex degree $n_H.$ For large $n_E$ we have that $\lambda_E \approx 2n_H^{-1/2},$ see \cite{hlw}.

{\em (vii)} Let $G_m$ be the zig-zag product between $E$ and $H,$ see \cite{rvw} for the definition and basic properties of the zig-zag product. The vertex set of $G_m$ is $V(E) \times V(H),$ it is a $d_H^2$--regular graph, and the absolute value of its second eigenvalue satisfy $\lambda_{G_m} \leq \lambda_E + \lambda_H.$

{\em (iix)} By construction $G_m$ has $n_S\cdot m^d$ vertices and there is a bijection $\phi_1:V(E) \rightarrow \{0,1,\ldots,n_S-1\} \times \{0,1,\ldots, m-1\}^d.$ We can assume that $V(E)$ is represented by $\{0,1,\ldots, n_E-1\}$ and do this fast by modulo arithmetic. Compose $\phi_1$ with $\phi_2 : (h,x_1,\ldots,x_d) \mapsto ms_h+(x_1,\ldots,x_d).$ This gives a bijection to the translated points of $S$ situated inside big cubes whose middle points are inside $Q.$ These points are the small red points inside the jagged green triangle in Figure~\ref{fig:B}.
In Figure~\ref{fig:B} we had $m=2.$ The local transition to from $m=2$ to $m=6$ is drawn in Figure~\ref{fig:A}.

{\em (ix)} Each of the points of $(S+\mathbf{Z}^d) \cap mP$ are given by a vertex of $G_m,$ and we refer to them as the {\em relevant} vertices. The other vertices are {\em irrelevant} vertices. Note that as $m \rightarrow \infty$ the proportion of irrelevant vertices tends to a constant 
$\varepsilon >0.$

\bpoint{Running the MCMC}

{\em (x)} We random walk on $G_m$ and for each vertex we test if it's relevant or not by evaluating linear inequalities of $mP$. The irrelevant ones are discarded, and for large $m$ this is a constant proportion. We may choose $H$ to achieve that $\lambda_H+2n_H^{-1/2} < 0.99$ for large enough $m.$ This shows that the second eigenvalue of $G_m$ is bounded away from 1 by a constant as $m \rightarrow \infty,$ or equivalently:

\tpoint{Theorem} The constructed Markov chain is rapidly mixing.

\bpoint{General remarks}
If the proportion of discarded vertices is too high, it can be pushed arbitrarily close to zero by replacing $P$ by $m'P$ for some big positive integer $m'.$ In practical examples this is a bad idea; it is surprisingly efficient in examples from algebraic statistics to discard a huge proportion, like 99.99\%, because the Markov basis inside the main construction might be much easier to calculate then. If we are prepare to discard even more, the blue triangle in Figure~\ref{fig:B} might be replaced by a box or some other polytope containing $P$ that is easy to understand -- compare to the extension methodology in optimisation theory.

An obvious extension is to not scale $P$ by $m$ in all axes but to have different scalings depending on the direction. It works fine to modify the main construction to achieve this, but we didn't spell it out since it doesn't add anything conceptually.

It is tempting to ask if the construction could be recursed, with $H$ coming from it, and repeat several times. This is possible, but not necessary. The best constructions of expanders already contains a repeated use of the zig-zag product to build them from simple starting graphs. Also referring to the theory of expanders, see \cite{hlw} and \cite{rvw}, one shouldn't think of building them in memory but rather as building an algorithm that is used in each instance of the random walk taking a step. Some of these algoritms are very explicit and we believe that for particular classes of problems from algebraic statistics it should be possible to write down the explicit zig-zag products with the explicit $H$ graphs from Markov basis, and put this into industry grade software as R or Matlab.

\bpoint{Remarks related to algebraic statistics}\label{algi}
The following interpretation is conceptually correct, but not completely correct. Consider the problem of $4 \times 4$ contingency tables with row and column sums $m$ Let the {\em floor} of a table be that we replace each entry or $r$ by $\lfloor 3r/m \rfloor.$ The entries of the floor of a table will be $0,1,2$ or $3.$ Note that we will never have $2,2,2,2$ in a row or column of the floor of a table. The graph $H$ in the construction should be thought of as a Markov chain on the floors of the contingency tables. The expander $E$ is on $m^{3^2}$ vertices because there are $3^2$ degrees of freedom. To get a table back from the floor and the expander, one would multiply the floor by $m,$ add the expander values on the top $3 \times 3$ square, and then equate to get the correct row and column sums. Sometimes numbers would become negative when equating, and this corresponds to falling outside the polytope, that is, being an irrelevant vertex.

To be completely clear, we should also clarify that the theory of expansion focus on when the asymptotic distribution is uniform. Sometimes in algebraic statistics that is not the original setting, and it's not clear how to transfer results to that setup. For the hypergeometric setting there is a natural way to extend the weights to the irrelevant vertices, and in (small) computational experiments by the author it seems like that it's mixing better than in the uniform case. On occasions the author have also conjectured that to hold in general.

In \ref{nm} we explained why the original Markov bases setup doesn't provide rapid mixing. In practice it can anyways be interesting to understand the mixing time, in particular in the hypergeometric setting. For a transition matrix $T$ with stationary distibution $v$ numerical evidence of slow mixing is provided by a vector $w$ that is far from $v$ while $Tw$ is close to $w.$ In this remark we explain an heuristic how to find such $w.$ The Markov bases are frequently derived from Gr\"obner bases, see \cite{dss}. In a Gr\"obner basis there is an order on the variables, and that order can be specified as an weight order on the monomials. The vertices of our graphs have coordinates, and the coordinates are exponents of the monomials. Thus, there is an induced weight $\omega$ on the coordinates, and it is actually a linear form taking rational values. For some $r$ about half of the vertices gets a weight $\omega$ larger than $r.$ In many problems of algebraic statistics there are symmetries that makes the choice of $r$ easy. Define $w$ to have value $1$ for vertices with weight $\omega > r$ and value $-1$ for vertices with weight $\omega < r.$
Rescale $w$ if necessary. In \ref{nm} we cut a polytope containing the graph by an codimension one hyperplane --- this is an explicit construction of that where $\omega = 0$ defines the hyperplane.

\end{document}